\documentclass[12pt]{article}
\usepackage{amssymb}
\usepackage{amsmath,amsthm}

\begin{document}

\title{\textbf{\ Turbulence in Fluid and Plasma: Search for a New Paradigm}}
\author{Dhurjati Prasad Datta\thanks{Corresponding author, email: dp\_datta@yahoo.com} \\ 
 Department of Mathematics \\
  University of North Bengal,
Siliguri, West Bengal, India, Pin 734013}
\date{}
\maketitle

\baselineskip=15.5pt

\begin{abstract}
A first principle explanation of the origin of intermittency and nonlinear structure formation in the Lagrangian velocity increments of a turbulent flow is presented in the context of a scale invariant analytical formalism that is being developed recently. The copious generation of power laws and nonlinear exponents in the structure functions are shown to follow quite naturally in the present formalism.
\end{abstract}

\begin{center} 

\end{center}

\begin{center} 
\end{center}


\section{Introduction}

Understanding the emergence and  subsequent control of turbulent flows in fluid motion described by the Navier-Stokes equation is of considerable interest recently. Starting from the seminal contributions of Richardson and Kolmogorov, there have been great advances in various directions \cite{frisch}, for instance, the statistical explanations of energy casceding and dissipation, emergence of intermittency pioneered by Kolmogorov \cite{kol} and Obkubov cite{obh} vis a vis the analysis of dynamic equations of hydrodynamical forces given by the Navier Stokes equation based on Shell models  \cite{shell}, Stochastic Langevin equations \cite{pinton} and so on. An interesting area of investigation is the exploration of the relationship of turbulence and chaos in deterministic systems. Although the question of existence of chaos in Navier Stokes equation is still open, there have been some recent advances in the existence proofs of simpler equations such the Euler equation \cite{li}. 

Origin of intermittency in the dynamics of turbulence is a very active field of current research \cite {pinton,bacry,vul,ben}. Such studies are recently being addressed mainly mainly in the framework of the Lagrangian approach: experimentally, one follows the motion of a single tracer particle and the relevant quantity of study is the increments in time of its velocity fluctuations $\Delta_{\tau} v(t)=v(t+\tau) -v(t)$, where $v(t)=v(x_0, t)$ parametrises the velocity of the fluid particle initially at $x_0$ along a Lagrangian trajectory \cite{pinton}. It is suggested that the anomalous scaling in the Lagrangian velocity increments could be traced back to the existence of long time correlations  in the particle accelearations that relate to the hydrodynamic forces  driving the particle motion in the turbulent flow. These long time correlations are expected to be of fundamental importance in 
understanding both the dynamical origin of turbulence and the associated statistical modelling.

A key issue in the intermittency research, in the context of a turbulent flow, is to look for a {\em dynamical explanation} of the anomalous scaling exponents in the structure functions of the velocity increments $<|\Delta_{\tau} v(t)|^q>\sim \tau^{\xi(q)}$ because of the extreme non-gaussianity of the corresponding probability distributions. A linear exponent of the form $\xi(q)=q/2$ corresponds to the Browninan motion and nonintermittent Kolmogorov flow. A nonlinear behaviour of $\xi(q)$ reveals the presence of multifractality leading to intermittency \cite{bacry}. Thus intermittency is linked to the nature and distribution of small scale {\em structure formations} on the fluid particle trajectory leading to intense bursty activity over a more quiet background. 

In this paper, we present a novel scale invariant dynamical (more generally, analytical) framework which can support such {\em nonlinear structure} formation leading to {\em localised bursty activity} in a {\em very natural} way \cite{sd1, sd2, ad, dp1, dp2}. We also identify the scaling exponents with the Hausdorff dimensions of a class of Cantor like fractal sets those are shown to arise spontaneously on the trajectory of a Lagrangian particle along the turbulent flow. {\em This approach is likely to offer a general framework for addressing key issues related to chaotic and turbulent motions in geometrically and intuitively much simpler manner compared to those available in the literature.} It should also offer {\em new insights} into the emergence and evolution of multifractal structures in various fields such as turbulence, finance, web traffic and so on.

To highlight the actual status of the present work against the current activities \cite {pinton,bacry,vul,ben} in the areas of intermittency and multifractality in fully developed turbulence, let us mention that understanding the emergence and dynamics of intermittent multifractal structures in turbulent flows are considered to be an open issue. Current researches mainly aim at developing analytic concepts under various special assumptions\cite{bacry,vul} and computational techniques \cite{pinton} in analysing  relevant questions and interpreting experimental facts. This work now presents a {\em novel analytic mechanism} for generating copiously intermittent and multifractal structures by exploiting a new class of {\em nonsmooth} solutions to linear differential equations (Sec.2) that becomes available when the possibility of a new level of {\em inversion induced nonlinear increments} for a dynamical variable is formulated rigorously (Sec.3). The relevant concepts and results are  formulated and developed in the context of a  Langevin stochastic differential equation when the noise term is assumed to have a special scale invariant form (Sec.2). The random fluctuations in the stochastic noise term is then related to the inversion mediated jump distributions on a Cantor like fractal set {\em leading to nontrivial velocity increments}. A one parameter class of such Cantor sets is shown to arise naturally in the present scale invariant formalism applied to the turbulent flow. The bursty intermittent behaviour is then interpreted as due to occurrence of Devil's staircase functions as ultrametric {\em weights (valuations)} assigned naturally to such multifractal Cantor sets. As a consequence, this extended scale  invariant   analysis and the associated dynamical framework turns out  to offer a natural framework to study dynamics of trubulent flows (Sec.4) and also of other complex systems. In short, this paper offers a new approach towards ushering in a paradigm shift in addressing key questions related to turbulent flows. The standard classical analytic framework might very well be extended to a scale invariant formalism that has been shown here to offer more natural interpretations to observed facts such as intermittent structure formations in turbulence. The present work is, however, mainly of theoretical interest and does not yet make any contact to experimental data. 

\section{Motivation}
The formulation of the {\em scale invariant analysis} is motivated by the possible existence  of a class of {\em nonsmooth} solutions of the simplest  ode of the form \cite{dp2}
\begin{equation}\label{se}
{\frac{dv}{dt}}=-{\frac{1}{T_0}}v
\end{equation}

The formalism of the scale invariant analysis is developed originally in the context of an extended real analysis incorporating dynamically active infinitesimals living in a non-archimedean space \cite{ad}. However, we outline here {\em essential features of the said analysis using the language of turbulence}. As a consequence,  $v(t)=v({\bf r}_0, t)$ here is assumed to represent, for definiteness, the Lagrangian velocity field of a tracer fluid particle in a turbulent flow. In ordinary sense, the above equation of course represents a {\em laminar phase} that relaxes in time $T_0$. 

The origin of nonsmooth  solutions can now be explained as follows. Let us begin by noting that a solution is {\em nonsmooth} in the sense that it could at most be realised as an approximate solution in the framework of ordinary calculus and also in the classical mechanics. However, in an extended {\em non-archimedean, ultrametric space}, such a solution can indeed be accommodated as a smooth solution. Our aim here is to highlight how such nonsmooth solutions could play a significant role in generating intermittent behaviour in a turbulent flow.

Notice that in ordinary calculus,  an independent real variable $t$, considered here as time, is supposed to {\em  undergo changes (increments)} only by linear shifts of the form $t\rightarrow t+\Delta t$, and eq(\ref{se}) is known to have the unique solution $v(t)=v(0) e^{-{\frac{t}{T_0}}}$. Suppose that a localised disturbance  is created (by stirring, say,) on the initial laminar flow of fluid, say,  at time $T_0$. The subsequent turbulent flow is modelled phenomenologically as a stochastic Langevin equation of the form 

\begin{equation}\label{sfe1}
{\frac{dv}{dt}}=-{\frac{1}{T_0}}v +{\frac{dF}{dt}}
\end{equation}
where forcing noise term is represented as $dF={v}{\frac{d\log \tau}{dt}}dt$. Such a stochastic modelling is expected to describe the fluctuating motion of the tracer particle in a turbulent flow \cite{pinton} with a very low viscosity and hence at a very high Reynold number $Re= l_0v_0/\nu$, where $l_0$ is the integral scale, $v_0$ stands for the rms (root mean square) velocity fluctuation at the integral scale and $\nu$ is the kinematic viscosity. For a white noise, the above equation leads to the Brownian motion.  However, in the present approach, the basic Gaussian random variate is assumed instead to be the {\em  scale invariant (stochastic) time} $\tau$,   and so the noise is log-normal and hence is coloured. The scale invariant time $\tau$ has a value of order O(1) when the ordinary Newtonian time $t$ is of O(1). However, $\tau$ must approach to infinity as $t\rightarrow \infty$. The above ansatz for the noise term now could be interpreted as revealing the existence of long time correlations in the time variable $t$ itself due to the induced fluctuating variable $\tau$, thus having the {\em novel} representation $t\sim T_0(1+\eta+\log \tau)$, $\eta$ being a growing (O(1)),  ordinary (real) variable, and $T_0\rightarrow \infty$. Notice that $t$ reduces to the ordinary time variable viz, $t\sim T_0(1+\eta)$ when $T_0\sim O(1)$. This long time correlated fluctuation then, in turn, will be inherited by the underlying velocity fluctuations as well. As a consequence, {\em  noise term is incorporated into the very definition (ie, the intrinsic structure) of time itself} (cf. next section for more precise justifications). This is one of the {\em salient features} of the formalism. Removing the laminar part $v_0=v(0)e^{-\frac{t}{T_0}}, \ t\sim T_0(1+\eta)$, the renormalised velocity fluctuation $\tilde v=vv_0^{-1}$ is now given formally by the scale invariant equation 

\begin{equation}\label{sfe}
d\tilde v=\tilde v d\log \tau
\end{equation}

Eq(\ref{sfe}), although should be interpreted as a (coloured) stochastic equation, however, in the following is treated  as an ODE (deterministic system) written {\em instead over a positive measure Cantor set} : a fattened real number system $\cal R$, each point of which is replaced by a zero measure {\em infinitesimal} Cantor set $C_p$ belonging to a one parameter $p$  family of such sets. The elements of the Cantor set $C_{p0}$, that replaces the ordinary 0 of $R$,  are new {\em infinitesimally} small elements in $\cal R$, known as {\em soft zeros} having dynamic properties. The scale invariant ODE (\ref{sfe}) is well defined over the gaps, of all possible sizes, of such a  Cantor set \cite{ad}. {\em The randomness in the fluctuations is therefore transferred to the dynamical formation of Cantor like sets in the back ground of a turbulent flow.} 

Now to formulate motion over a Cantor set rigorously, we make an {\em assumption} that the scale invariant time $\tau$ undergos increments by {\em random inversions} of the form $\tau\rightarrow \tau^{-1}$, in the neighbourhood of $\tau=1$ (say) . Let $\tau_{-}= 1 - \eta(t)$ and $\tau_{+}= 1 + \tilde\eta(t)$ such that the {\em inversion rule} $\tau_-\tau_+=1$, or more generally, $\tau_+\tau_-^{1+a}=1$, $a$ being a real parameter indicating the size of the increment (jump), is satisfied. The inversions are realised randomly in the sense that a fluid particle at a site $\tau_-$ jumps over instantneously to the site $\tau_+$ when $\tau_-$ and $\tau_+$ are separated by a gap of size $a$.  Because of small scale (Cantor set like) structures, the exact specification of $\tau_+$ in the fattened real number set $\cal R$ is also not precisely known. The inversion rule then defines a small scale multiplicative structure close to  $\tau=1$, that arises dynamically as the ordinary time $t$ approaches asymptotically the limiting value  $\infty$ (cf, Sec 3 for details). The multiplicative scaling exponent $a$ in $\tau_+=\tau_-^{-1-a}$ is an ultrametric valuation defined by 

\begin{equation}\label{value1}
a(\tau)=\underset{\tau_-\rightarrow 0}{\lim} [{\frac{\log \tau_+ - \log \tau_-^{-1}}{\log \tau_-^{-1}}}] =\underset{\tau_-\rightarrow 0}{\lim} [{\frac{\log \tau_+/ \tau_-^{-1}}{\log \tau_- ^{-1}}}]
\end{equation}
as  $t\rightarrow \infty$. As a consequence, and also because of scale invariance, the ordinary connected line segment neighbourhood of $T_0$ is now extended to a positive measure Cantor like set ${\tilde C}_p$, realised as an ultrametric space with the absolute value $a$. {\em Transition between granular like points of such a Cantor set, realised as end points of the gap,  are mediated by seemingly stochastic inversion induced jumps} (of sizes determined by the associated gap sizes) as defined above, so that the corresponding multiplicative jump differential measure is defined by $\Delta_j \tau =\log a(t), \ t\rightarrow \infty$. Notice further that the valuation $a$ is a {\em locally constant function} defined by $d a(t)=0$ which is a direct consequence of the fact that $d \log\log \tau_+=d\log\log \tau_-^{-1}$, so that $a(t)$ has a globally intermittent behaviour given by a Devil's stairecase function (i.e. a Cantor function) indicating {\em significant fluctuations} in the neighbourhood of a point of the ultrametric Cantor set (ie. the set over which the valuation undergoes actual variations). However, $a(t)$, although it is  constant over the scale of $t$, ie, ${\frac{da}{dt}}=0$, its intrinsic variability is revealed in the smaller $\log\log t$ scale, ie, 

\begin{equation}\label{vari}
\log t^{-1}\frac{d\tau }{d\log t^{-1}}=\tau.
\end{equation}

To summarise the {\em main advantage} of our analytic approach over others is the dynamic generation of infinite numbers of cascaded  smaller scales of the form $\lambda_1{\frac{l}{L}}<l_n<\lambda_2{\frac{l}{L}}$ for every pair of bounded intervals of the form  $0<\lambda_1<\lambda_2<<1$ in a scale invariant manner and thereby assuring formation of a host of Cantor like sets, all of which are distributed following certain multifractal measures.

\section{Formalism}
\subsection{Revisiting Ordinary Calculus}

Let us recall that analysis of motion, from simple and/or laminar motion to more complex chaotic and/or turbulent flows, ultimately rests  on the basic concepts of differential calculus. Let us begin first by showing how a slight variation of conventional treatment and meaning of the statement $\underset{t \rightarrow 0} {\lim} \ t=0, \ t$ being a real variable can lead to potentially new results which would be significant for justifying the subsequent  scale invariant formalism. The ordinary meaning of the above limiting statement is, of course, that the real variable $t$ approaches linearly and uniformly with constant rate 1 to the limiting value 0 remaining always on the connected line segment $[0,t]$. 

This continuous and connected flow may, however, be broken by introducing a scaling of the form $f(t)=at, \ 0<a<1$ and subsequently iterating  infinite number of times. The original length (Lebesgue measure)  $t$ now reduces to 0 successively by jumps as $a^nt, $ as $n\rightarrow \infty$. Conventionally, one feels that these are only two possibilities with unique end result. 

Surprisingly, however, there  exists still another nontrivial possibility. To every $t$ and $a$ satisfying above constraints, one is free to choose an arbitrary finite $p: 1<p<a^{-1}$ such that the continuous linear measure of the flow $t\rightarrow 0$ is preserved at every level $n$, following a {\em measure conservation principle}

\begin{equation}\label{consv}
t=\{ta^{n\log t^{-1}}\}\times p^{n\log T}
\end{equation}
The scale invariance of the above equality is obvious. The linear flow (in the lhs), as it were, is transformed into a scale invariant nonlinear flow in logarithmic scales. As a consequence, the reduction of the linear measure 1, at each scale by the factor $a^{n\log t^{-1}}$ is corrected multiplicatively by creating a space in a verticle direction containing the fattened variable $T>1$ in the logarithmic scale and as a consequence, in the limit $n\rightarrow \infty$, and hence as $t\rightarrow 0$, the linear variable $t$ is transformed into a nonlinear flow pattern characterised by the scaling law $T^{-1}=t^{\frac{1}{s}}$, where $s=\frac{\log p}{\log a^{-1}}$ is the fractal (Hausdroff) dimension of a Cantor like set that may be assumed to have been formed spontaneously in the neighbourhood of 0, because of the limit process. The precise dynamical setting of realising a cascade of scales of the form $a^n$ and the associated (daughter) multiplicative scales $p^n$ will be presented in the next section. We note that, because of the inversion rule, the original variable $t$ on its approach towards 0, would get reversed (ie, experience a bounce) close to 0, so that the corresponding scale invariant infinitely large fattened variable $T_s= \tilde T^{-1}= tT^{-1} $ will live and vary on a Cantor set of fractal dimension $\tilde s=\frac{\log (1/ap)}{\log p}$. This, in turn, would have a very important consequence: {\em the value 0 is actually never reached; the linear flow along the positive $t$ axis, in the interval $(0,t), \ t\downarrow$, will be transferred in a scale invariant manner, to a new branch in the inverted variable $T_s\sim t^{\tilde s}$, by inversion induced jumps. The subsequent (linear) flow of $T_s$ toward 0 will then again experience a bounce, and so on leading to a very intricate multifractal tree like structure.} The multifractality follows from the existence of a one parameter family of multiplicative scale factors $p(a)=pa$ depending on $a$.

 To have a more concrete appreciation of the above analysis, let us suppose that points  of the line segments [0,1] have the structure  of tiny infinitesimally small beads aligned horizontally. Under an application of the contraction the length of the original segment is reduced to $a$. How do the infinitesimal beads align on the contracted segment? Since no part of the original segment is removed, the beads would climb one above other to create a space dynamically; and this process would continue in the limit $n\rightarrow \infty$, and ultimately one realises {\em a fractured fractal like structure} in the neighbourhood of 0, but in the vertical direction.
 
 In the following section, we briefly present the formal structure of the scale invariant analysis accommodating the above novel structure in a rigorous manner.
 
\subsection{Scale Invariant Analysis}

We state the basic definitions introducing scale invariant dynamic infinitesimals, and the associated non-archimedean norms and also outline briefly the natural emergence of intermittent structures induced by Cantor's stairecase functions (for details, refer to \cite{ad}) and their possible role in the asymptotic behaviours of a dynamic variable related to a turbulent flow.
 
By {\em dynamic infinitesimals} we mean a class of quantities living in a well defined space, which are arbitrarily small in comaprison to any finite element of the space, however, might have nontrivial O(1) effect on the long time evolutionary properties of such elements. For definiteness, the underlying space is assumed to be the set of real numbers $R$, and any finite real variable is denoted by $x\in R$. In a dynamical problem the variable $x$ is replaced by the ordinary time $t$. 

{\bf Definition 1:} Let $x\in I=[0,1]\subset R$ and $x$ be arbitrarily small, i.e., $x\neq 0$, but, nevertheless, $x\rightarrow 0^+$. Then there exists $\delta>0$ and a set of (positive) \emph{relative infinitesimals} $\tilde x$ relative to the scale $\delta$ satisfying $0<\tilde x<\delta\leq x $ and the \emph{inversion} rule $\tilde x/\delta\propto
\delta/x$. The associated scale invariant infinitesimals are defined by $
\tilde X=\underset{\delta\rightarrow 0^+}{\lim} \tilde x/\delta$.

As an example, for a continuous variable $x$ approaching $0^+$ and considered itself as a scale, a class of the relative infinitesimals are represented as $\tilde x \propto x^{1+l}(1+o(x)), \ 0<l<1$, corresponding to the real variable $x^{1-l}$, so that the corresponding scale invariant 
infinitesimals are defined by the formula $\tilde X=\lambda x^l +o(x^m), \ m>l$. Choosing the parameter $\lambda$ from an open set $U\subset (0,1)$, we get a class of relative infinitesimals $\tilde x$ belonging to an open subset of $(0, \delta)$, all of which are related by the inversion rule to  the real number $x^{1-l}$ for given  $x$ and $l$. Moreover, allowing $\lambda$ to vary from a countable number of disjoint open intervals $U_r, \ r=1,2,\ldots$, one can introduce a countable, disjoint class ${\tilde I}_r$ of relative infinitesimals, all of which are related by the inversion rules $\tilde x_r/\delta = \lambda_r\delta/x$, where $\lambda_r\in U_r$. Notice that $\underset{r}{\bigcup}{\tilde I}_r \subset (0,\delta)$. It  follows therefore that the infinitesimals actually reside in the gaps of an arbirarily chosen Cantor subset of the open interval $(0,\delta)$. As a consequence, the nonlinearity induced via the inversion rule in the apparently linear limiting motion of $x$ naturally leads to a large family of Cantor like fractal sets.  

Notice that in the limit $\delta \rightarrow 0$, the set of relative infinitesimals apparently reduces to the singleton $\{0\}$. However, the set of scale
invariant infinitesimals may  nevertheless be nontrivial. The set of relative infinitesimals, denoted $\mathbf{0}$, has the {\em asymptotic
representation} $\mathbf{0}=\{0, \delta \tilde X_r\}$, as $ \delta\rightarrow 0^+$.
For definiteness, the ordinary zero (0) is called the \emph{stiff} zero, when
the nontrivial infinitesimals are called \emph{soft or dynamic} zeros. The ordinary
real line $R$ is then extended over $\mathbf{R}=\{\mathbf{\tilde x}: \mathbf{
\tilde x}=x+\mathbf{0}, \ x\in R\}$, which as a field extension, and because of the Frobenius theorem, must be an infinite dimensional nonarchimedean space. 

{\bf Definition 2: Non-archimedean Norm} The mapping $v:\mathbf{0}\rightarrow I^+=[0,1]$
\begin{equation}\label{norm}
\phi(\tilde x)\equiv ||\tilde x|| := \underset{\delta\rightarrow 0^+}{\lim}%
\log_{\delta^{-1}} {\tilde X}^{-1}, \ \tilde x=\delta \tilde X \ \in \mathbf{%
0}
\end{equation}
\noindent together with $\phi(0)=0$ defines  a non-archimedean (ultrametric ) norm on the set of infinitesimals $\mathbf{0}$. Infinitesimals weighted with above absolute value are called
\emph{valued} infinitesimals. Clearly, the ordinary Euclidean norm for a dynamic infinitesimal would vanish trivially because $|\tilde x_r|=\delta|\tilde X_r|=0$, in the limit $\delta\rightarrow 0$.

Let us recall that an absolute value is non-archimedean if the triangle inequality is replaced by the stronger ultrametric inequality $v(a+b)\leq {\rm max}\{v(a), \ v(b)\}$. 
Ultrametricity now tells that the dynamic infinitesimals are indeed awardednontrivial {\em values} of the form $\phi(\tilde x_r)= a_rp^{-rs}$ for $p=p(\delta)>1$ and $a_r>0$ being an ascending sequence continuously on the closure of $\underset{r}{\bigcup}{\tilde I}_r$, so that the valuation $\phi(\tilde x)$ has the structure of an intermittent Cantor function associated to a Cantor set with fractal dimension $s$ that has formed dynamically in the infinitesimal neighbourhood of 0. {\em Since the ultrametric parameter $p$ generates a one parameter family of scale factors corresponding to a family of Cantor sets, the nonlinear structure that emerges out from the above analysis  of course corresponds to multifractality.}

One of the most important features of dynamic infinitesimals (by inversion, infinitely large scales) is that {\em these can actually lead to significant influences in the long time asymptotic properties of an evolving dynamical system}. As a result, a laminar flow can indeed be transformed into a complex turbulent flow over a sufficiently long time scale and also possibly in the presence of a slight disturbance.

To see the precise manner how the small scale nonlinear structure inhabited by the dynamic infinitesimals can induce an O(1) variation in the long time properties of any system, let us first note that the dynamic infinitesimals introduces a sort of a deformation into the standard structure of finite real numbers $x$ through the nontrivial valuation in the form

\begin{equation}\label{defm}
x\rightarrow \ X_\pm= x\cdot x^{\mp \phi(\tilde x)} 
\end{equation}
where the valuation $\phi$, realised as a Cantor's stairecase function, satisfies eq(\ref{vari}). The scale invariant factor $\tilde X_\pm=X_\pm/x$ corresponding to the {\em fattened} variable $X$, by assumption, undergo transition by inversions, so that the scale invariant increments in $\tilde X$ may be defined as  

\begin{equation}\label{inc}
\delta_s \tilde X=\log (X(x)/x)= \phi(\tilde x)\log x^{-1}
\end{equation} 
which matches with the linear increment $\delta_l \tilde X=|\tilde X-1|= \phi(\tilde x)\log x^{-1}$ when $x\rightarrow 0$. The corresponding linear increment for the fattened variable has the form $\delta_l X= |X-x|=\phi(\tilde x)x\log x^{-1}$ as $x\rightarrow 0$. Notice that {\em these are new level of increments which arise out of inversion mediated variations, that in the present formulation are activated in presence of dynamic infinitesimals leading to small scale structure formation which in turn would create long time correlations and very large intermittent fluctuations} \cite{pinton}.

\section{Applications to Turbulence}

To apply the above analytical results to  a turbulent flow, let us notice that the standard multifractality assumption suggests formation of cascades of eddies between the integral and dissipation scales converging towards a family of fractal sets $C_h$ with fractal Hausdorff dimension $D(h)$ for each scaling exponent $h$ from a bounded interval of the form  $(h_{\rm min}, h_{\rm max})$. The corresponding increments of the Lagrangian velocity fluctuations $\delta _{\tau} v= v(t+\tau)-v(t)$ are {\em assumed} to follow the scaling law $\delta _{\tau}v\sim v_0(\tau/\tau_0)^h$ for each scaling exponent $h$ from the bounded interval.  Here, $v_0$ and $\tau_0$ correspond to the integral scale characteristic rms velocity fluctuation and the corresponding integral time scale respectively. Notice that the derivation of such power laws is impossible in the ordinary (classical) framework of Navier-Stokes equation. Clearly, the scale invariant analytical formalism now offers us with a mechanism to {\em realise} this multifractal scaling {\em dynamically} in a most natural manner.

Notice that one can identify the scale invariant velocity increment $\delta_{\tau}v/v_0$ in the turbulent flow with the increment of the corresponding scale invariant fattened velocity $\delta_s \tilde V$ following the incremental law (\ref{inc}) in the limit of sufficiently small (stochastic) time increments $0<\tau<<\tau_0$, so that one obtains $\delta_{\tau} \tilde V\sim \phi(\tilde \tau)\log \tau^{-1}$. Now, recalling that the valuation $\phi$ enjoys the structure of a Cantor Stairecase function, it follows that $(\delta _{\tau}v)/v_0:= \delta_{\tau} \tilde V\sim \phi(\tilde \tau)\log \tau^{-1}$ scales as $(\delta _{\tau}v)/v_0\sim (\tau/\tau_0)^{h-\epsilon}$ for an $\epsilon>0$ (that arises form the logarithmic factor), and $h$ being the Hausdorff dimension of the Cantor set $\tilde C_h$ that emerges dynamically in association with the dynamic infinitesimals and the scale invariant limiting processes in the present formalism. The above scaling also reveals that the velocity undergoes variations only over such family of Cantor sets. Let us emphasize that {\em this power law velocity increment follows  not from  any hypothesis , but presented here as first principle  derivation from the scale invariant approach.} Notice that the set $\tilde C_h$ is the set that arises in association to the velocity field fluctuations, when $C_h$ is the set that is realised in physical space of cascaded eddies in the fluid medium.

Adopting the well accepted analysis available in literature \cite {frisch, pinton, bacry}, one then immediately concludes that the structure functions of the velocity fluctuations $S_q(\tau) =E(|\delta v|^q)\sim (\tau/\tau_0)^{\xi(q)}$, where the scaling exponent is a nonlinear concave function of $q$, and can be defined following Ref.\cite{frisch} as a Legendre transform of $D(h)$: $\xi(q)= \underset{h}{\rm inf}[qh+1-D(h)]$. Here, $D(h)>0$ is the fractal dimension of the set $C_h$ on which cascades of eddies accumulate when the indefinite interations are realised internally as $t\rightarrow \infty$. Notice, also that $1-D(h)$ actually stands for the {\em codimension} of the limit set $C_h$ in the 1 dimensional time axis.  The derivation of the scaling exponent $\xi(q)$ in the context of multifractal random walk models \cite{pinton, bacry} will be considered elsewhere.

To summarise, we have outlined a first principle derivation of nonlinear multifractal scaling exponents $\xi(q)$ of Lagrangian velocity fluctuations. The analysis is based on a simple reformulation of the definition of a limit process exposing how small scale structures are realised in the neighbourhood of a disturbance. When the laminar flow is stirred vigorously, the velocity $v(t)=v({\bf r_0}, t)$ of a fluid particle at $\bf r_0$ at time $t$ would have a spectrum of  fluctuating values $v(t+\tau)$ for scale invariant $0<\tau<<1$ living in a Cantor like set $\tilde C_h$ with fractal dimension $h$. Instead of relaxing quickly to the laminar phase, these fluctuations will exist as  stable structures corresponding to the scaling exponents $\xi(q)$. Conversely, the energy injected at the integral scale by stirring cascades through smaller scales following nonlinear convective terms of the Navier Stokes equation. These spatial fluctuations then induces Lagrangian velocity fluctuations in the time domain which can be interpreted (following the present analysis) as generated from intrinsic,  scale invariant stochastic variations in the time $t$ itself. As a consequence, multifractal eddy formations can be considered as physical menifestations  of {\em dynamical multifractals} those are shown to exist in the infinitesimal time scales. More detailed analysis and interpreting $\xi(q)$ in the light of experimental data will be considered separately.

Finally, we note that this scale invariant analytical formalism might also be crucial in elucidating the nature of the plasma turbulence observed in the space plasma and high temperature fusion plasma devices where the turbulence is anomalous. A detailed analysis of this  aspect will be reported in a separate publication.

\end{document}